\documentclass[12pt,a4paper]{article}

\usepackage{amsmath,amssymb}

\def\oM{{\overline{\cal M}}}

\def\Cbar{{\overline C}}
\def\Ihat{{\widehat I}}

\def\Ptilde{{\widetilde P}}
\def\cC{{\cal C}}
\def\oC{{\overline{\cal C}}}
\def\qed{{\hfill $\diamondsuit$}}
\def\cL{{\cal L}}
\def\CP{{{\mathbb C}{\rm P}}}
\def\Aut{{\rm Aut}}
\def\Z{{\mathbb Z}}

\def\C{{\mathbb C}}
\def\Q{{\mathbb Q}}

\def\sing{{\rm sing}}
\def\basic{{\rm basic}}

\usepackage{theorem}

\newtheorem{theorem}{Theorem}
\newtheorem{proposition}{Proposition}[section]

\newtheorem{lemma}[proposition]{Lemma}

{\theorembodyfont{\rmfamily}
\newtheorem{definition}[proposition]{Definition}
\newtheorem{example}[proposition]{Example}
\newtheorem{remark}[proposition]{Remark}
\newtheorem{notation}[proposition]{Notation}
}

\usepackage{graphicx}

\title{Universal cohomological expressions for singularities in
families of genus~0 stable maps}

\author{Maxim Kazarian\thanks{Steklov Mathematical Institute RAS,
National Research University Higher School of Economics},
Sergey Lando\thanks{National Research University Higher School of Economics,
Independent University of Moscow},
Dimitri Zvonkine\thanks{
Institut math{\'e}matique de Jussieu,
CNRS, 4 place Jussieu,
75005 Paris, France. E-mail: dimitri.zvonkine@gmail.fr. \newline
Partly supported by the ANR project ``Geometry and
Integrability in Mathematical Physics'' ANR-05-BLAN-0029-01.
}}

\date{}

\begin{document}

\maketitle

\begin{abstract}
We consider families of curve-to-curve maps that have no singularities except those of genus~0 stable maps and that satisfy a versality condition at each singularity. We provide a universal expression for the cohomology class Poincar\'e dual to the locus of any given singularity. Our expressions hold for any family of curve-to-curve maps satisfying the above properties.
\end{abstract}

\section{Introduction}

\subsection{A family of curve-to-curve maps}

A family of curve-to-curve maps is a triple of smooth complex manifolds $X$, $Y$ and $B$ and a commutative diagram
\begin{equation} \label{Eq:maps}
\begin{picture}(50,25)(0,25)
\put(0,35){$X$}
\put(45,35){$Y$}
\put(22,0){$B$}
\put(14,38){\vector(1,0){25}}
\put(9,30){\vector(2,-3){11}}
\put(45,30){\vector(-2,-3){11}}
\put(22,42){$f$}
\put(5,19){$p$}
\put(42,19){$q$}
\end{picture}
\end{equation}

\vspace{15\unitlength}
\noindent
satisfying the following properties.

\begin{enumerate}
\item The map $p$ is a flat family of reduced not necessarily compact nodal curves.
\item The map $q$ is a flat family of reduced not necessarily compact smooth curves.
\item The map~$f$ is proper.
\end{enumerate}

In Section~\ref{Ssec:sing} we will introduce a list of allowed singularities for the map~$f$. These will be precisely the singularities of genus~0 stable maps. For each singularity we will also require a versal deformation property: informally, the map~$f$ should locally present all possible deformations of each singularity. Assuming that~$f$ has allowed singularities only and satisfies the versal deformation property, we will compute the cohomology class of every singularity locus.

The forthcoming paper~\cite{KLZ2} extends the results of the present one to the study of multisingularities.

\subsection{Thom polynomials: a motivation}

Let $f:X \to Y$ be a sufficiently generic map between
two complex compact manifolds. Denote by
$$
c(f) = c(f^*(T_*Y))/c(T_*X) = 1 + c_1(f) + c_2(f) + \dots,
$$
the full Chern class of the map~$f$, $c_i(f)\in H^{2i}(X)$.

Given a singularity type $\alpha$, one can consider the
locus of points $x \in X$ such that $f$ has a singularity of
type $\alpha$ at~$x$. This locus is a submanifold, and we will
denote the Poincar\'e dual cohomology class of its closure by~$[\alpha] \in H^*(X,\Z)$.

In~\cite{Thom} Thom proved that {\em for every singularity
type $\alpha$ and for every generic map~$f$ we have $[\alpha] = P_\alpha(c_1(f), c_2(f), \dots)$,
where $P_\alpha(c_1, c_2, \dots)$ is a polynomial independent of~$f$.}

For instance, we have $[A_1] = c_1(f)$, thus $P_{A_1} = c_1$.

The polynomials $P_\alpha$ are called {\em Thom polynomials}.
It is, in general, quite hard to compute them (Thom's proof is
not constructive); however computing even a single
Thom polynomial may provide a solution for many seemingly unrelated
problems of enumerative geometry.

This paper grew out of an attempt to ``apply'' Thom's theorem
to the map~$f$ of Eq.~\eqref{Eq:maps}, a case to which it is
not applicable since the genericity condition for~$f$ brakes down due to existence of nonisolated singularities.

In Section~\ref{Ssec:sing} we describe all
singularity types that can appear in a family of genus~0 stable maps.
In Section~\ref{Ssec:classes} we introduce a family of cohomology
classes (larger than just the classes $c_i(f)$) that we will
call {\em basic classes}. Then we prove that the class
determined by every singularity type has a universal polynomial
expression in terms of the basic classes, independent of the family~$f$.
Moreover, we give a simple way to compute these expressions.

\subsection{Local models}
Consider two families of curve-to-curve maps.
\begin{equation} \label{Eq:maps}
\begin{picture}(50,25)(0,25)
\put(0,35){$X$}
\put(45,35){$Y$}
\put(22,0){$B$}
\put(14,38){\vector(1,0){25}}
\put(9,30){\vector(2,-3){11}}
\put(45,30){\vector(-2,-3){11}}
\put(22,42){$f$}
\put(5,19){$p$}
\put(42,19){$q$}
\end{picture}
\qquad \qquad
\begin{picture}(50,25)(0,25)
\put(-1,35){$X'$}
\put(44,35){$Y'$}
\put(21,0){$B'$}
\put(14,38){\vector(1,0){25}}
\put(10,30){\vector(2,-3){11}}
\put(45,30){\vector(-2,-3){11}}
\put(22,42){$f'$}
\put(2,19){$p'$}
\put(42,19){$q'$}
\end{picture}
\end{equation}

\vspace{2.5em}

Let $y \in Y$ and $b = q(y) \in B$; $y' \in Y'$ and $b' = q'(y') \in B'$. Let $X_0$ be a closed subset of $f^{-1}(y)$ and $X'_0$ a closed subset of $(f')^{-1}(y')$. In practice, $X_0$ will be either a point or a connected component of~$f^{-1}(y)$ contracted by~$f$.

\begin{center}
\includegraphics[width=10em]{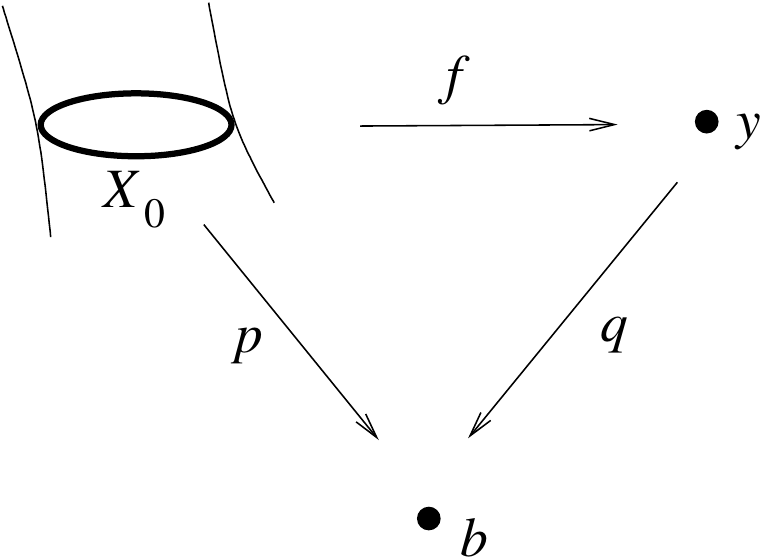}
\end{center}

\begin{definition}
We say that the family $(X',Y',B')$ is a {\em local model} for the family $(X,Y,B)$ if there exist three holomorphic submersions $\varphi_X:X \to X'$, $\varphi_Y:Y \to Y'$, $\varphi_B:B \to B'$ defined in a neighborhood of $X_0$, $y$, and $b$ respectively, identifying $X_0$ with $X_0'$ and commuting with the maps $p,q,f$.
\end{definition}

\subsection{Singularity types}
\label{Ssec:sing}

Let $x \in X$ be any point, $b = p(x) \in B$ and
$y  = f(x) \in Y$ its images in $B$ and~$Y$, respectively,
and $X_b$ the fiber of~$p$ through~$x$. We will describe several types of singularities of $f$ at~$x$, namely, ramification points, nodes, and contracted components. These are the only singularities that can arise
in families of genus~$0$ stable maps. A map $C\to\CP^1$ of a genus~$0$ nodal curve~$C$ to~$\CP^1$
is said to be {\it stable\/} if its automorphism group is finite. In other words, a map is stable
if each irreducible component of the domain~$C$ that is contracted to a point has at least three
points of intersection with other irreducible components of~$C$.

For each singularity type we will define a versality condition using local models. Both points and nodes are isolated singularities, while
contracted components are nonisolated ones. Universal polynomials for isolated singularities in the absence
of non-isolated ones were found in~\cite{KazLan0,KazLan}. They determine the singularity classes modulo
cohomology classes supported on the subvariety in~$X$ consisting of contracted components.
Including contracted components completes the computation of universal polynomials for genus~$0$ curves.

\subsubsection{Ramification points}

\begin{definition} \label{Def:singA}
We say that $f$ presents a {\em singularity of
type $A_k$ at~$x$} if $x$ is a smooth point of $X_b$ and $f|_{X_b}$ has a ramification of order
$k+1$ at~$x$.
\end{definition}

The local model for the $A_k$ singularity is given by the family
$$
v (u, \gamma_1, \dots, \gamma_{k-1}) = u^{k+1} + \gamma_1 u^{k-1} + \dots + \gamma_{k-1} u.
$$
Here:
\\ $(\gamma_1, \dots, \gamma_{k-1})$ is a system of coordinates in~$B'$, \\ $(\gamma_1, \dots, \gamma_{k-1}, u)$ a system of coordinates in~$X'$,
\\ $(\gamma_1, \dots, \gamma_{k-1}, v)$ a system of coordinates in~$Y'$.

\begin{definition} \label{Def:versal1}
Let $x \in X$ be a point at which $f$ presents an $A_k$ singularity. We say that the family $(X,Y,B)$ is {\em versal} at~$x$ if it is locally modeled at the neighborhood of~$x$ by the family $(X',Y',B')$ above.
\end{definition}

\subsubsection{Nodes}
\begin{definition}
We say that $f$ presents a {\em singularity of
type $I_{k_1,k_2}$ at $x$} if $x$ is a node of $X_b$ and $f$ has ramifications of orders $k_1$ and $k_2$ on the branches of $X_b$ at~$x$.
\end{definition}

A local model for the singularity $I_{k_1, k_2}$ is given by the family
$$
\beta = u_1 u_2,
$$
$$
v = u_1^{k_1} + \gamma^{(1)}_1 u_1^{k_1-1} + \dots + \gamma^{(1)}_{k_1-1} u_1 + u_2^{k_2} + \gamma^{(2)}_1 u_2^{k_2-1} + \dots + \gamma^{(2)}_{k_2-1} u_2.
$$
Here:
\\ $(\gamma^{(1)}_1, \dots, \gamma^{(1)}_{k_1-1}, \gamma^{(2)}_1, \dots, \gamma^{(2)}_{k_2-1}, \beta)$ is a system of coordinates in~$B'$,
\\ $(\gamma^{(1)}_1, \dots, \gamma^{(1)}_{k_1-1}, \gamma^{(2)}_1, \dots, \gamma^{(2)}_{k_2-1}, u_1, u_2)$ a system of coordinates in~$X'$,
\\ $(\gamma^{(1)}_1, \dots, \gamma^{(1)}_{k_1-1}, \gamma^{(2)}_1, \dots, \gamma^{(2)}_{k_2-1}, \beta, v)$ a system of coordinates in~$Y'$.

\begin{definition} \label{Def:versal2}
Let $x \in X$ be a point at which $f$ presents an $I_{k_1,k_2}$ singularity. We say that the family $(X,Y,B)$ is {\em versal} at~$x$ if it is locally modeled at the neighborhood of~$x$ by the family $(X',Y',B')$ above.
\end{definition}

\subsubsection{Contracted components}
Stable maps of genus~0 can also present non-isolated singularities.
Suppose that $x$ happens to lie on an irreducible component of $X_b$ contracted by~$f$. Then we denote by $\cC$
the {\em connected} component of $f_b^{-1}(y) \subset X_b$ containing $x$. Note that it necessarily includes the irreducible component of $X_b$ containing~$x$, but may include other irreducible components as well. We call $\cC$ the {\em contracted part} of~$x$. Finally we call {\em branches} the connected components of the neighborhood of $\cC$ in $X_b \setminus \cC$. If, on $\cC$, we mark its intersection points $x_1, \dots, x_\ell$ with the branches, it becomes a stable curve.
\begin{figure}
\begin{center}
\begin{picture}(150,157)(0,-5)
\put(3,4){\includegraphics{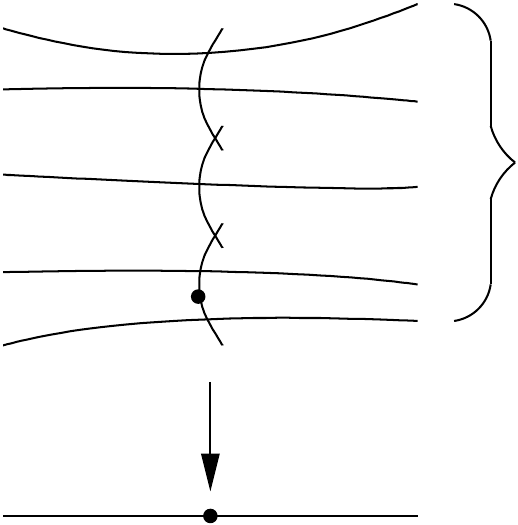}}
\put(50,30){$f_b$}
\put(58,-7){$y$}
\put(54,52){$k_1$}
\put(50,81){$k_2$}
\put(51,149){$k_\ell$}
\put(50,66){$x$}
\put(70,115){$\cC$}
\put(160,107){branches}
\end{picture}
\end{center}
\caption{The contracted part~$\cC$ of a point~$x$. Branches are indicated, together
with the orders $k_1,\dots,k_\ell$ of the restriction of the function~$f_b$ to the
corresponding branch}
\end{figure}
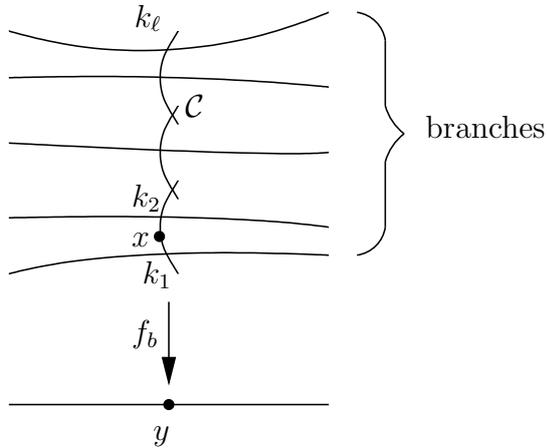


\begin{definition} \label{Def:I}
For $\ell \geq 3$, we say that $f$ presents a {\em singularity of
type $I_{k_1,k_2, \dots, k_\ell}$ at $x$} if
\begin{itemize}
\item
$f$ is constant on the irreducible component of $X_b$ containing~$x$;
\item
the contracted part $\cC$ has genus~0;
\item
there are $\ell$ branches and they can be numbered from 1 to~$\ell$ in such a way that $f$ has ramification orders $k_1, \dots, k_\ell$ at their intersection points with~$\cC$.
\end{itemize}
\end{definition}


To define a local model for this singularity we need to introduce the space of relative stable maps; however, the local model may actually be only a subspace of this space, see Definition~\ref{Def:versal3}.

Let $\oM_{0;k_1, \dots, k_\ell}(\CP^1;\infty)$ be the space of stable maps of genus~0 to $\CP^1$, defined up to an additive constant, and relative to $\infty$ with ramification profile $(k_1, \dots, k_\ell)$. This space has a distinguished sub-orbifold $\oM^Z_{0; k_1, \dots, k_\ell}(\CP^1;\infty)$ that we will call the {\em zero locus}. It parametrizes maps from curves with $\ell+1$ components: one contracted rational component and $\ell$ more rational components that intersect the contracted one and on which the map equals $z^{k_i}$, $1 \leq i \leq \ell$. As an orbifold, the zero locus is a gerb with base $\oM_{0;\ell}$ and group $\prod_{i=1}^\ell \Z / k_i \Z$. Let

$p’:\oC_{0;k_1, \dots, k_\ell}(\CP^1;\infty) \to \oM_{0;k_1, \dots, k_\ell}(\CP^1;\infty)$ be the universal curve,

$q’: \CP^1 \times  \oM_{0;k_1, \dots, k_\ell}(\CP^1;\infty) \to  \oM_{0;k_1, \dots, k_\ell}(\CP^1;\infty)$ the universal target curve, and

$f’ : \oC_{0;k_1, \dots, k_\ell}(\CP^1;\infty) \to \CP^1 \times \oM_{0;k_1, \dots, k_\ell}(\CP^1;\infty)$ the universal map.

\noindent
The union of contracted parts of fibers over the zero locus $\oM^Z_{0;k_1, \dots, k_\ell}(\CP^1;\infty)$ forms the zero locus $\oC^Z_{0;k_1, \dots, k_\ell}(\CP^1;\infty)$ in the universal curve.

Let $x \in X$ be a point at which $f$ presents an $I_{k_1, \dots, k_\ell}$ singularity and $\cC$ the corresponding contracted component. Let $b \in \oM^Z_{0; k_1, \dots, k_\ell}(\CP^1;\infty)$ be the point of the zero locus such that the contracted component is isomorphic to~$\cC$. Consider a local chart of $\oM_{0;k_1, \dots, k_\ell}(\CP^1;\infty)$ at~$b$. Now, let $B'$ be any sub-manifold of this local chart that is transversal to the zero locus and of dimension at least
$$
\dim \oM_{0;k_1, \dots, k_\ell}(\CP^1;\infty) - \dim \oM^Z_{0;k_1, \dots, k_\ell}(\CP^1;\infty)
=
\sum k_i.
$$
Let $X'$ be a neighborhood of $\cC$ in the preimage $(p')^{-1}(B')$ and $Y'$ a neighborhood of $y’ = f’(x’)$ in the preimage $(q')^{-1}(B')$.

\begin{definition} \label{Def:versal3}
Let $x \in X$ be a point at which $f$ presents an $I_{k_1, \dots, k_\ell}$ singularity and $\cC$ the corresponding contracted component. We say that $f$ is {\em versal} at~$x$ if there exists a family $(X',Y',B')$ as above such that $(X,Y,B)$ is locally modeled on $(X',Y',B')$ at a neighborhood of $\cC$.
\end{definition}

\subsection{Versal families} \label{Ssec:Examples}

\begin{definition}
A {\em genus 0 versal family of maps} is a family of maps that only has singularities of types $A_k$ and $I_{k_1, \dots, k_\ell}$ for $\ell \geq 2$ and satisfies the versality condition for each singularity.
\end{definition}

Genus~0 versal families are precisely the object of our study.
Before proceeding we give several examples of genus~0 versal families.

\begin{example}
Consider the space $\oM_{0;d}(\CP^1)$ of stable genus~0 degree~$d$ maps to $\CP^1$ without marked points. The universal map over this space is a genus~0 versal family.
\end{example}

\begin{example}
Given a line bundle $L \to B$ over a smooth base~$B$, consider the space
$$
(L \setminus \{\mbox{zero section} \}) \times_{\C^*} \oM_{0;d}(\CP^1) \to B.
$$
Here $\C^*$ acts on $\oM_{0;d}(\CP^1)$ via the natural action on the target $\CP^1$. The universal map over this space is a genus~0 versal family.
\end{example}

\begin{example}
Given a principal ${\rm PSL}(2,\C)$ bundle $G \to B$ over a smooth base~$B$, consider the space
$$
G \times_{{\rm PSL}(2,\C)} \oM_{0;d}(\CP^1) \to B.
$$
Here ${\rm PSL}(2,\C)$ acts on $\oM_{0;d}(\CP^1)$ via the natural action on the target $\CP^1$. The universal map over this space is a genus~0 versal family.
\end{example}

\begin{example}
Consider the space $\oM_{0;k_1, \dots, k_\ell}(\CP^1;\infty)$ of stable genus~0 maps to $\CP^1$ relative to $\infty$ with ramification profile $(k_1, \dots, k_\ell)$. Consider the universal map over this space. Remove the $\infty$ section in the target curve and its preimages in the source curve. The family thus obtained is a genus~0 versal family.
\end{example}

\begin{example}
Given a line bundle $L \to B$ over a smooth base~$B$, consider the space
$$
(L \setminus \{\mbox{zero section} \}) \times_{\C^*} \oM_{0;k_1, \dots, k_\ell}(\CP^1;\infty) \to B.
$$
Here $\C^*$ acts on $\oM_{0;k_1, \dots, k_\ell}(\CP^1;\infty)$ via the natural action on the target $\CP^1$. Consider the universal map over this space. Remove the $\infty$ section in the target curve and its preimages in the source curve. The family thus obtained is a genus~0 versal family.
\end{example}

\begin{example}
Consider the space $\oM_{g;d}(C)$ of degree~$d$ genus~$g$ stable maps to a given smooth curve~$C$. Consider the universal map over this space. Remove from the universal map all the source curves that contain contracted components of genus greater than 0. Also remove the images of these curves in the target. The family thus obtained is a genus~0 versal family.
\end{example}

\begin{example}
Let $z=z_1, \dots, z=z_\ell$ be a given set of points on~$\CP^1$ with affine coordinate~$z$. Consider the family of maps
$$
f(z) = \sum_{i=1}^\ell  \left[
\left(\frac{u_i}{z-z_i} \right)^{k_i} + \gamma_{i,1} \left(\frac{u_i}{z-z_i} \right)^{k_i-1} + \cdots + \gamma_{i,k_i-1} \left(\frac{u_i}{z-z_i} \right) \right].
$$
This family is extended to $u_i=0$ in the following way. The source $\CP^1$ acquires a bubble with global coordinate $w_i$. The point $w_i=0$ is attached to $z=z_i$. The function on the bubble is given by
$$
w_i^{k_i} + \gamma_{i,1} w_i^{k_i-1} + \cdots + \gamma_{i,k_i-1} w_i.
$$
This family is a genus~0 versal family.
\end{example}

\begin{example}
This is an example of a family that is NOT versal. Let $z'=z'_1, \dots, z'=z'_{\ell_1+1}$ be a given set of points on the projective line with global coordinate $z'$ and  $z''=z''_1, \dots,z''= z''_{\ell_2+1}$ be a given set of points on another projective line with global coordinate $z''$. The point $z'_{\ell_1+1}$ is identified with the point $z''_{\ell_2+1}$.
Consider the following family of maps on the nodal curve thus obtained:
$$
f(z') = \sum_{i=1}^{\ell_1}  \left[
\left(\frac{u'_i}{z'-z'_i} \right)^{k'_i} + \gamma_{i,1} \left(\frac{u'_i}{z'-z'_i} \right)^{k'_i-1} + \cdots + \gamma_{i,k'_i-1} \left(\frac{u'_i}{z'-z'_i} \right) \right] + C';
$$
$$
f(z'') = \sum_{i=1}^{\ell_2}  \left[
\left(\frac{u''_i}{z''-z''_i} \right)^{k''_i} + \gamma_{i,1} \left(\frac{u''_i}{z''-z''_i} \right)^{k''_i-1} + \cdots + \gamma_{i,k''_i-1} \left(\frac{u''_i}{z''-z''_i} \right) \right] + C'';
$$
where the constants $C'$ and $C''$ are such that $f(z'_{\ell_1+1}) = f(z''_{\ell_2+1})=0$. This family is extended to $u'_i=0$ and $u''_i=0$ as in the previous example.

This family is not versal. Indeed, even though it satisfies the versal property for the singularity $I_{k'_1, \dots, k'_{\ell_1}, k''_1, \dots, k''_{\ell_2}}$, a generic map of this family also has an $I_{1,1}$ singularity that does not satisfy the versal property because the node is never smoothened.
\end{example}

\subsection{Singularity loci} \label{Ssec:SingLoci}
Every singularity from the list in Sec.~\ref{Ssec:sing} determines
a cycle in $X$, namely the closure of the set of points
$x \in X$ where $f$ presents a given singularity.
By abuse of notation we will usually denote in the same
way the singularity type and the corresponding cycle in $X$.

For instance, the cycle $I_{1,1}$ is the set of nodes of
the fibers $X_b$.

The versal deformation property insures that the codimension of the $A_k$ locus equals $k$, that of the $I_{k_1, k_2}$ locus equals $k_1 + k_2$, and that of the $I_{k_1, \dots, k_\ell}$ locus equals $\sum k_i$.

In the notation $I_{k_1, \dots, k_\ell}$
the indices form a multiset, that is, the order of the indices
is immaterial. Since the cycle $I_{k_1, \dots, k_\ell}$
is the {\em closure} of the set where $f$ presents
the corresponding singularity, the ramification orders
of~$f$ at the attachment points of the branches can
actually be greater than (but not smaller than)~$k_i$.
Thus, for instance, the cycle $I_{1,2}$ presents a
self-intersection along the cycle $I_{2,2}$. Therefore
the normalization of $I_{1,2}$ contains two copies
of $I_{2,2}$.

Let $|\Aut\{k_1, \dots, k_\ell \}|$ be the number of permutations $\sigma \in S_\ell$ such that $k_{\sigma(i)} = k_i$ for every~$i$.

Fix an order of the indices $k_1, \dots, k_\ell$.
Let $\Ihat_{k_1, \dots, k_\ell}$ be the space of couples
$(x \in I_{k_1, \dots, k_\ell}, \mbox{numbering of the branches})$,
such that the ramification order at the branch number~$i$
is equal to~$k_i$ (greater than or equal to~$k_i$ in the closure).
Then  $\Ihat_{k_1, \dots, k_\ell}$ is
a finite nonramified covering of the normalization of
$I_{k_1, \dots, k_\ell}$. Its degree equals $|\Aut \{k_1, \dots, k_\ell \}|$.




\begin{notation}
We denote by $a_k \in H^{2k}(X,\Q)$ the cohomology class Poincar\'e dual to
the cycle $A_k$. We denote by $i_{k_1, \dots, k_\ell} \in H^{2 \sum k_i} (X,\Q)$
the cohomology class Poincar\'e dual to the cycle
$I_{k_1, \dots, k_\ell}$ multiplied by $|\Aut\{k_1, \dots, k_\ell \}|$.
\end{notation}


\begin{remark}
Including the factor $|\Aut\{k_1, \dots, k_\ell \}|$ in the definition of cohomology classes will greatly
simplify the formulas in the sequel.
\end{remark}

\begin{notation} \label{Not:uv}
We denote by $U \colon \Ihat_{k_1, \dots, k_\ell} \to I_{k_1, \dots, k_\ell}$ the covering map and by $V \colon \Ihat_{k_1, \dots, k_\ell} \to \oM_{0, \ell+1}$ the natural map from the covering $\Ihat$ to the moduli space of curves:
$$
I_{k_1, \dots, k_\ell} \stackrel{U}{\longleftarrow}
\Ihat_{k_1,\dots, k_\ell} \stackrel{V}{\longrightarrow} \oM_{0,\ell+1}.
$$
Given a class $\alpha \in H^*(\oM_{0,\ell+1})$,
we denote by $\alpha i_{k_1, \dots, k_\ell}$ the class
$U_*V^*(\alpha)$.
\end{notation}

Note that the notation is coherent: if $\alpha = 1$ we have
$U_* V^*(1) = i_{k_1, \dots, k_\ell}$ as defined previously,
because the map~$U$ has degree $|\Aut \{k_1, \dots, k_\ell \}|$.

\begin{definition} \label{Def:Sing}
The cohomology classes $a_k$ and $\alpha i_{k_1, \dots, k_\ell}$ in
$H^*(X,\Q)$ are called {\em singularity classes}.
\end{definition}

\subsection{Basic classes}
\label{Ssec:classes}

Now we are going to introduce a different set of cohomology classes in~$X$.
We will call them {\em basic classes}.

\begin{definition}\label{Def:classes} \

Let $\cL_X \to X$ be the line bundle whose fiber
over $x \in X$ is the {\em cotangent line to $X_b$ at~$x$}.
This line bundle is well-defined outside $I_{1,1}$ and
can be uniquely extended to $I_{1,1}$ (in the standard way).
{\em We denote by $\psi = c_1(\cL_X) \in H^2(X,\Q)$
its first Chern class}.

Let $\cL_Y \to X$ be the line bundle whose fiber
over $x \in X$ is the {\em cotangent line to $f(X_b)$ at~$y=f(x)$}.
This line bundle is well-defined everywhere
and {\em we denote by $\xi = c_1(\cL_Y) \in H^2(X,\Q)$ its first Chern class}.
\end{definition}

Consider the cycle $I_{1,\dots,1}$ with $\ell$ subscripts.
Over its $\ell!$-sheeted covering
$\Ihat_{1,\dots,1}$ consider the $\ell$ line bundles
whose fibers are the {\em cotangent lines to the branches
at their intersection points with the contracted part}.
Denote their first Chern classes by $\nu_1, \dots, \nu_\ell$. In a more traditional notation, we have $\nu_i=V^*\psi_i$, $i=1\dots,\ell$.

Let $\alpha \in H^*(\oM_{0,\ell+1})$ be a cohomology class.
We denote by $\alpha \delta_{m_1, \dots, m_\ell}$ the cohomology class
$U_* (\nu_1^{m_1} \cdots \nu_\ell^{m_\ell} V^*\alpha)$.

\begin{definition} \label{Def:Basic}
The classes $\psi^k$, and $\alpha \delta_{k_1, \dots, k_\ell}$ are called
{\em basic cohomology classes}.
\end{definition}


Our aim is now to express singularity classes in terms
of the basic classes and vice versa.

\subsection{Main results}
\label{Ssec:results}

We start with an explicit equality.

\begin{theorem} \label{Thm:induction}
For every $m \geq 1$ we have
$\prod\limits_{r=1}^m (r \psi -\xi) = a_m + P_m$,
where $P_m$ is a linear combination of terms of the form
$\xi^q \psi^p i_{k_1, \dots, k_\ell}$. Moreover $P_1 = 0$ and
$$
P_m = \left(
\sum_{\ell \geq 2} \frac1{\ell!} \sum_{k_1 + \dots + k_\ell = m} \!\!\!\!
k_1 \cdots k_\ell \; i_{k_1, \dots, k_\ell} \right)
\; + \;
(m \psi-\xi) \, P_{m-1}
$$
for $m \geq 2$.
\end{theorem}

Sample calculations for $1 \leq m \leq 6$
are given in the appendix.

The equality of Theorem~\ref{Thm:induction} does not express
a singularity class via basic classes or a basic class
via singularity classes. However it is the
key result in the proof of the following theorem.

\begin{theorem} \label{Thm:SingBas}
Every singularity class can be expressed as a linear
combination of basic classes multiplied by powers of~$\xi$.
Every basic class can be expressed as a linear
combination of singularity classes multiplied by powers of~$\xi$.
\end{theorem}

The proof of this theorem is constructive, that is, the expressions can be effectively computed. In the appendix we write out the expressions for all classes up to codimension~5.

Finally, in the following theorem, we give an explicit formula for certain coefficients in the expressions of basic classes in terms of singularity classes.

Introduce the following polynomials in variables $x_k$:
$$
X_m = \sum_\ell \frac1{\ell!}
\sum_{
\substack{k_1, \dots, k_\ell\\
\sum k_i = m-\ell+2
}}
\frac{m!}{(m-\ell+2)!} \prod_{i=1}^\ell k_i x_i.
$$

Denote by $\alpha_\ell \in H^{\rm top}(\oM_{0;\ell+1})$ the
cohomology class Poincar\'e dual to a point.

\begin{theorem} \label{Thm:CompCycles}
Choose $m$ and $k_1, \dots, k_\ell$ so that
$m+2 = \ell+ \sum\limits_{i=1}^\ell k_i$.
The coefficient of $\alpha_\ell i_{k_1, \dots, k_\ell}$ in
the expression for $\psi^m$ is equal to the coefficient
of $x_{k_1} \cdots x_{k_\ell}$ in $X_m$, that is, to
$$
\frac1{|\Aut \{ k_1, \dots, k_\ell \}|}
\frac{m!}{(m-\ell+2)!} \prod k_i.
$$

Choose $m_1, \dots, m_s$ and $k_1, \dots, k_\ell$ so that
$2s+\sum\limits_{j=1}^s m_j  = \ell+ \sum\limits_{i=1}^\ell k_i$.
The coefficient of $\alpha_\ell i_{k_1, \dots, k_\ell}$
in the expression for $\alpha_s \delta_{m_1, \dots, m_s}$
is equal to the coefficient of $x_{k_1} \cdots x_{k_\ell}$
in $X_{m_1} \cdots X_{m_s}$.
\end{theorem}

In Section~\ref{Sec:CompCyc} we relate these coefficients
with Okounkov and Pandharipande's {\em completed cycles}.

\section{Computations with singularity and basic classes}

\subsection{Proof of Theorem~\ref{Thm:induction}}

The proof goes by induction. The equality $\psi-\xi = A_1$
follows from the fact that $df$ is a section of the line
bundle $\cL_X \otimes \cL_Y^*$ (see Definition~\ref{Def:classes})
that has a simple zero precisely on the $A_1$ locus.

Assume that the equality is true for $m-1$ and prove it
for~$m$. By the induction assumption we have
$$
\prod_{r=1}^m (r \psi - \xi)
= (m \psi - \xi) (a_{m-1} + P_{m-1}).
$$
The term $(m \psi - \xi) P_{m-1}$ is a part of the
expression for $P_m$. Thus our main task is to
compute the product $(m\psi - \xi) a_{m-1}$.

To do that, note that $d^mf$ is a section of the line
bundle $\cL_X^{\otimes m} \otimes \cL_Y^*$ restricted
to the locus $A_{m-1}$. Thus we must describe the components
of the zero locus of $d^mf$ and their multiplicities.

The obvious component of the zero locus is
$A_m \subset A_{m-1}$. On this component the
section $d^mf$ has a simple zero as follows from the versality condition for $A_m$. Indeed, in the local model of Definition~\ref{Def:versal1} the stratum $A_{m-1}$ is parametrized by
one variable~$a$ via
$$
u^{m+1} + \gamma_1 u^{m-1} + \dots \gamma_{m-1} u = (u-a)^m (u+ma)+ma^m;
$$
the section $d^m f|_{u=a}$ in this parametrization is just~$(m+1)a$, so it has a simple zero at $a=0$.

The other components are the cycles $I_{k_1, \dots, k_\ell}$
that lie in the closure of $A_{m-1}$ and are of codimension~$1$
in $A_{m-1}$. We are going to show that all cycles $I_{k_1, \dots, k_\ell}$ of codimension~$m$, that is, $\sum k_i = m$, lie in the closure of $A_{m-1}$ and that the vanishing multiplicity of $d^m f$ on the cycle $I_{k_1, \dots, k_\ell}$ equals $k_1 \cdots k_\ell$. To do that, we study the local models for each of these singularities.

We will assume that $\ell \geq 3$ leaving the simpler and very similar case $\ell = 2$ to the reader
(see also~\cite{KazLan}).

Let $\oM_{0;k_1, \dots, k_\ell}(\CP^1;\infty)$ be the space of stable maps of genus~0 to $\CP^1$, defined up to an additive constant, and relative to $\infty$ with ramification profile $(k_1, \dots, k_\ell)$. Denote by $\oC_{0;k_1, \dots, k_\ell}(\CP^1;\infty)$ the universal curve and let $\oC^Z_{0;k_1, \dots, k_\ell}(\CP^1;\infty) \subset \oC_{0;k_1, \dots, k_\ell}(\CP^1;\infty)$ be the zero locus in the universal curve (see the discussion before Definition~\ref{Def:versal3}).

\begin{lemma}
The natural section $d^mf:A_{m-1}\to\cL_X^{\otimes m} \otimes \cL_Y^*$  on the stratum $A_{m-1} \subset \oC_{0;k_1, \dots, k_\ell}(\CP^1;\infty)$ has a vanishing of order $k_1 \cdots k_\ell$ along $\oC^Z_{0;k_1, \dots, k_\ell}(\CP^1;\infty)$, whenever $\sum k_i =m$.
\end{lemma}

From the lemma we conclude that
$$
(m \psi - \xi) a_{m-1} = a_m + \sum_\ell \sum_{
\substack{k_1 \leq \dots \leq k_\ell\\
\sum k_i = m}
}
\!\!\!
k_1 \cdots k_\ell \; [I_{k_1, \dots, k_\ell}]
$$
$$
= a_m + \sum_\ell \frac1{\ell!}\sum_{
\substack{k_1, \dots, k_\ell\\
\sum k_i = m}
}
\!\!\!k_1 \cdots k_\ell \; i_{k_1, \dots, k_\ell}.
$$

This holds for any family that satisfies the versality condition of Definition~\ref{Def:versal3}. Theorem~\ref{Thm:induction} is proved. \qed

\paragraph{Proof of the lemma.} We say that a point of $\oC_{0;k_1, \dots, k_\ell}(\CP^1;\infty)$ is
{\em generic} if the contracted part is smooth. We will study the vanishing order at a generic point.

Introduce a coordinate~$z$ on the contracted part
and denote by $x,z_1, \dots, z_\ell$ the marked point
and the intersection points with the branches. It is proved in~\cite{ELSV} that $\oM_{0;k_1, \dots, k_\ell}(\CP^1;\infty)$ is a cone over its zero locus and that the fiber of this cone is parametrized by coordinates $u_i$, $1 \leq i \leq \ell$ and $a_{ij}$, $1 \leq j \leq k_i-1$ if we write the stable map in the following form:
\begin{equation} \label{Eq:HurwitzCoords}
f(z) = \sum_{i=1}^\ell \left[
\left(\frac{u_i}{z-z_i}\right)^{k_i} +
a_{i,1}\left(\frac{u_i}{z-z_i}\right) ^{k_i-1} +
\dots +  a_{i,k_i-1} \left(\frac{u_i}{z-z_i}\right)
\right].
\end{equation}

%
%

Let
$$
f(z) = \frac{(z-x)^m}{\prod_{i=1}^\ell (z-z_i)^{k_i}}.
$$
This is the unique, up to transformations $f \mapsto cf+b$, rational function with poles of orders
$k_i$ at the points $z_i$ and such that
$f'(x) = f''(x) = \dots = f^{(m-1)}(x) =0$.
Expand $f$ in the form~\eqref{Eq:HurwitzCoords} and denote by $\bar u_i$ and $\bar a_{ij}$ the coefficients thus obtained. Further, let
$$
K = \mbox{LCM}(k_1, \dots, k_\ell)
$$
and
$$
r_i = \frac{K}{k_i} \quad \mbox{ for } \quad 1 \leq i \leq \ell.
$$

Then the stratum $A_{m-1}$ in $\oM_{0;k_1, \dots, k_\ell}(\CP^1;\infty)$ consists of $d = k_1 \cdots k_\ell /L$ irreducible components parametrized by a complex parameter $c$ in the following way:
$$
u_i = \alpha_i c^{r_i} \bar u_i, \quad a_{ij} = \alpha_i^j c^{j r_i} \bar a_{ij}.
$$
Here $(\alpha_1, \dots, \alpha_\ell)$ is a collection of roots of unity, $\alpha_i^{k_i} =1$. We should choose one such collection in each of the $d$ orbits of $\Z/K \Z$ in $\Z/k_1 \Z \times \dots \times \Z/ k_\ell \Z$ and this choice determines the irreducible component of the stratum $A_{m-1}$.

In each component of the stratum, the $m$th derivative at~$x$ of the function corresponding to parameter~$c$ is equal to $c^K \cdot f^{(m)}(x)$. Thus the vanishing order of the $m$th derivative at $c=0$ equals $K$ for each component. Since there are $k_1 \cdots k_\ell /K$ components, we get the total vanishing order equal to $k_1 \cdots k_\ell$. \qed

\subsection{Proof of Theorem~\ref{Thm:SingBas}}

\subsubsection{Classes in $\oM_{0;\ell+1}$}

The definition of both basic and singularity classes
involves a class $\alpha \in H^*(\oM_{0;\ell+1})$. Before
proving the theorem we must introduce notation for
these classes. They will be described using marked trees.

\begin{definition} \label{Def:trees}
An {\em $\ell$-tree} or just a {\em tree} is a rooted tree with $\ell$ leaves;
the valencies of all vertices except the leaves and the
root are at least~3; the valency of the root is~1 (but we
don't call it a leaf). A {\em marked $\ell$-tree}
or just a {\em marked tree} is an $\ell$-tree whose
all vertices except the root are marked with nonnegative integers.
\end{definition}

The picture below shows a marked tree.

\begin{center}
\
\includegraphics{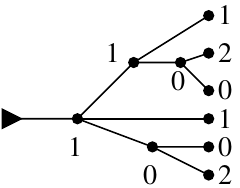}

\end{center}

If $T$ is a marked tree we will denote by $t$ the underlying
tree obtained by forgetting the markings.


Take a marked tree $T$ and number its leaves in an arbitrary way. This tree with numbered leaves determines a cohomology class in $\oM_{0;\ell+1}$ as follows. Consider the set of curves in $\oM_{0;\ell+1}$ whose dual graph\footnote{The dual graph of a curve is obtained
by replacing every irreducible component by a vertex,
every node by an edge, every marked point by a leaf,
and the point $x$ by the root.} is isomorphic to~$t$, the marked point $x$ being the root of the tree. Every interior vertex of~$t$ corresponds to an
irreducible component of the curve. On this irreducible component
there is exactly one marked point that leads to~$x$
(in other words, from every vertex~$v$ of~$t$ exactly
one edge leads to the root). The first Chern class of the
cotangent line bundle to this point is denoted by $\psi_v$. The cohomology  class that we assign to~$T$ with numbered leaves is given by
$[T] = \prod_{v\in V(T)} \psi_v^{m_v}$ supported on the boundary stratum corresponding to~$t$, where the product runs over all internal vertices~$v$ of~$t$ and~$m_v$ is the corresponding marking. Note that this class does not depend on the markings of the leaves.

Using the class $[T]$ we can now assign to a marked tree $T$ both a singularity class and a basic class. These classes will not depend on the numbering of the leaves of $T$, but will depend on their markings.
Denote by $m_1, \dots, m_\ell$ the integer markings on the leaves.

\paragraph{Basic class.}
The basic class that we assign to~$T$ is given by
$$
[T]_\basic = U_*\left(
V^* [T] \prod_l \nu_l^{m_l}
\right),
$$
where $U \colon \Ihat_{1, \dots, 1} \to I_{1, \dots,1}$ and $V \colon \Ihat_{1, \dots, 1} \to \oM_{0; \ell+1}$ are as in Definition~\ref{Def:Basic}.

\paragraph{The singularity class.}
The singularity class assigned to~$T$ is
$$
[T]_\sing = U_* (V^*[T]),
$$
where $U \colon \Ihat_{m_1+1, \dots, m_\ell+1} \to I_{m_1+1, \dots, m_\ell+1}$ and $V \colon \Ihat_{m_1+1, \dots, m_\ell+1} \to \oM_{0; \ell+1}$ are as in Notation~\ref{Not:uv} and Definition~\ref{Def:Sing}.

\begin{example}
If $T$ is the tree
\begin{center}
\
\includegraphics{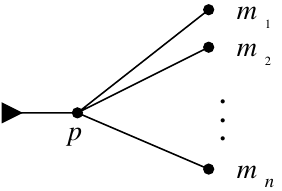}

\end{center}
then $[T]_\sing = \psi^p i_{m_1+1, \dots, m_\ell+1}$,
$[T]_\basic = \psi^p \delta_{m_1, \dots, m_\ell}$.
\end{example}

\begin{notation}
If $T$ is the tree
\begin{center}
\
\includegraphics{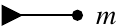}
\end{center}
we let, by convention, $[\includegraphics{stick.pdf}]_\basic = \psi^m$, $[\includegraphics{stick.pdf}]_\sing = a_m$,
$[\includegraphics{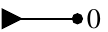}]_\basic = [\includegraphics{1sing.pdf}]_\sing  = 1$.
\end{notation}

Thus every singularity class and every basic class can
be represented as a linear combination of trees.

\begin{remark}
Keel~\cite{Keel} proved that the cohomology classes
Poincar\'e dual to the boundary strata of $\oM_{0;\ell+1}$
span the whole cohomology ring of $\oM_{0;\ell+1}$. Therefore
our tree notation is sufficient to express any class, but
is strongly redundant. For instance, whenever
the integer $d_v$ assigned to an interior vertex
of $T$ is greater than the valency of~$v$ minus 3, the corresponding
class vanishes. As another example of redundancy, we have
$[T_1]_\sing = 3 [T_2]_\sing$, where
\begin{center}
\
\includegraphics{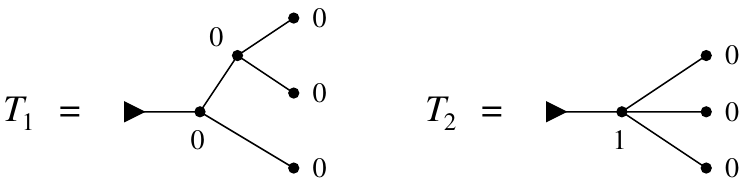}
\end{center}
However the classes in our expressions naturally appear
in the form of marked $\ell$-trees, and we chose not to simplify them further,
since we do not know any simple nonredundant expression
for the result of such simplification.
\end{remark}

We finish this section with a definition that we will need soon.

\begin{definition}
Suppose $\ell$ marked trees $T_1, \dots, T_\ell$ are assigned to
the leaves of a marked $\ell$-tree~$T$.
The {\em substitution} of $T_1, \dots, T_\ell$
into $T$ is the tree obtained by erasing the
leaves of $T$ and gluing the roots of $T_1, \dots, T_\ell$
into the vertices on which the leaves grew. Note that every
vertex of the substitution (except the root) inherits
an integer from either $T$ or one of the $T_i$'s, but not
both.
\begin{center}
\
\includegraphics{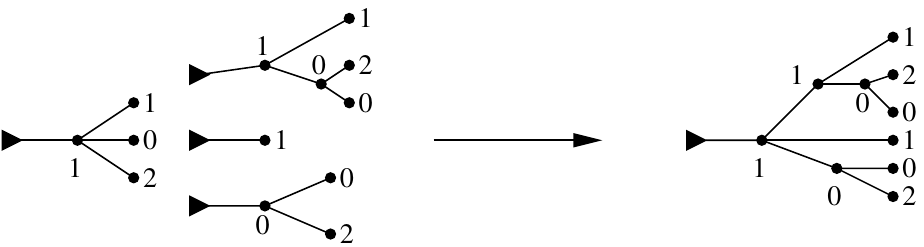}

\end{center}
\end{definition}

If instead of the trees $T_1, \dots, T_\ell$ we
have $\ell$ linear combinations of trees we can, of course, extend
the operation of substitution by multi-linearity.

\subsubsection{Basic classes via singularity classes}

Recall the list of basic classes: $\psi^k$ and $[T]_\basic$,
where $T$ is a marked $\ell$-tree, $\ell \geq 2$.
We must express every class like that as a linear
combination of singularity classes multiplied by powers of~$\xi$.

\paragraph{Expressions for $\psi^m$.} Theorem~\ref{Thm:induction}
gives an expression for $\prod\limits_{r=1}^m (r \psi - \xi)$ as
a linear combination of singularity classes multiplied by powers
of $\xi$. Writing
\begin{eqnarray*}
\psi &=& (\psi - \xi) + \xi,\\
\psi^2 &=&  \frac12 \cdot (2 \psi - \xi)(\psi- \xi)
+ \frac32 \xi \cdot (\psi - \xi)
+ \xi^2,\\
\psi^3 &=&  \frac16 \cdot (3 \psi - \xi)(2 \psi - \xi)(\psi- \xi)
+ \frac{11}{12} \xi \cdot (2 \psi - \xi)(\psi- \xi)
+ \frac74 \xi^2 \cdot (\psi - \xi) + \xi^3,
\end{eqnarray*}
and so on, we obtain similar expressions for powers of $\psi$.

\begin{example}
We have
$$
\psi = a_1 + \xi, \qquad
\psi^2 = \frac12 a_2 + \frac14 i_{1,1} +\frac32\xi a_1 + \xi^2.
$$
More computations are given in the appendix.
\end{example}

To sum up, for every $m$ there exists a linear combination
$L_m$ of trees with coefficients in $\Q[\xi]$ such that
$\psi^m = [L_m]_\sing$.

\paragraph{The general basic class.}
Consider class $[T]_\basic$ and denote by
$m_1, \dots, m_\ell$ the integers on the leaves of~$T$.

\begin{proposition} \label{Prop:substitution}
Let $\widehat T$ be the substitution of
$L_{m_1}, \dots, L_{m_\ell}$ into the tree~$T$.
Then the class $[T]_\basic$ is equal to $[{\widehat T}]_\sing$.
\end{proposition}

\paragraph{Proof.}
This proposition is almost obvious. It suffices to note
that every class $\nu_i$ plays the role of the class $\psi$
on the corresponding branch, and therefore we can use
the expressions for $\psi^m$ to express $\nu_i^{m_i}$
for each~$i$. \qed

\begin{example}
Let us compute the expression of $\delta_{0,1,2}$
in singularity classes. We have
$$
\delta_{0,1,2} = \left[ \begin{minipage}[c]{4.1em}
\includegraphics{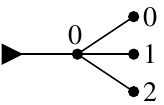}
\end{minipage} \right]_\basic,
$$
$$
1 =  \left[ \begin{minipage}[c]{2.7em}
\includegraphics{1sing.pdf}
\end{minipage} \right]_\sing,
\qquad
\psi = \left[ \begin{minipage}[c]{2.7em}
\includegraphics{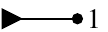}
\end{minipage} \right]_\sing
 + \xi
\left[ \begin{minipage}[c]{2.7em}
\includegraphics{1sing.pdf}
\end{minipage} \right]_\sing,
$$
$$
\psi^2 = \frac12
\left[ \begin{minipage}[c]{2.7em}
\includegraphics{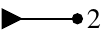}
\end{minipage} \right]_\sing
 + \frac14
\left[ \begin{minipage}[c]{4.1em}
\includegraphics{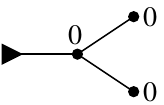}
\end{minipage} \right]_\sing
+ \frac32 \xi
\left[ \begin{minipage}[c]{2.7em}
\includegraphics{A1sing.pdf}
\end{minipage} \right]_\sing
+ \xi^2
\left[ \begin{minipage}[c]{2.7em}
\includegraphics{1sing.pdf}
\end{minipage} \right]_\sing
.
$$
Substituting the last three expressions into the first tree,
we obtain a linear combination of 8 trees that simplifies
to
$$
\delta_{0,1,2} = \frac12 i_{1,2,3} + \frac32 \xi i_{1,2,2}
+\frac12 \xi i_{1,1,3} + \frac52 \xi^2 i_{1,1,2}
+ \xi^3 i_{1,1,1} +
$$
$$
+
\frac14\left[
\begin{minipage}[c]{5em}
\includegraphics{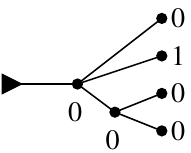}
\end{minipage}
\right]_\sing
+\frac14 \xi \left[
\begin{minipage}[c]{5em}
\includegraphics{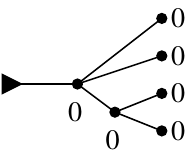}
\end{minipage}
\right]_\sing .
$$
\end{example}

\subsubsection{Singularity classes via basic classes}
Let $T$ be a marked $\ell$-tree with the markings of the leaves
equal to $m_1, \dots, m_\ell$. We call $m_1 + \cdots + m_\ell$
the {\em weight} of~$T$. It
follows from the previous section that the expression
of $[T]_\basic$ via singularity classes has the form
$$
[T]_\basic = \frac1{m_1! \cdots m_\ell!} [T]_\sing +
\mbox{ lower weight terms}.
$$
Therefore these expressions form a triangular change
of basis and can be inverted. Theorem~\ref{Thm:SingBas} is proved.

\bigskip

The expressions for the simplest singularity classes
are computed in the appendix.

\subsection{Proof of Theorem~\ref{Thm:CompCycles}}

\paragraph{The expression for $\psi^m$.}
The first claim of Theorem~\ref{Thm:CompCycles} gives
an explicit formula for the coefficient of
the class $\alpha_\ell i_{k_1, \dots, k_\ell}$ in the
expression of $\psi^m$, where $\alpha_\ell$ is the
class of a point in $\oM_{0;\ell+1}$. We will prove
this formula using Theorem~\ref{Thm:induction};
however, since the terms we are interested in
do not contain the class $\xi$, we can reduce the
formula of Theorem~\ref{Thm:induction} modulo~$\xi$.
We obtain:
$$
m! \psi^m = a_m + \Ptilde_m  \quad (\mbox{mod } \xi),
$$
where $\Ptilde_1 = 0$, and
$$
\Ptilde_m = \left(
\sum_{n \geq 2} \frac1{\ell!} \sum_{k_1 + \dots + k_\ell = m} \!\!\!\!
k_1 \cdots k_\ell \; i_{k_1, \dots, k_\ell} \right)
\; + \;
m \psi \, \Ptilde_{m-1}
$$
for $m \geq 2$.

Thus we see that the term $i_{k_1, \dots, k_\ell}$
first ``appears'' in $\Ptilde_{m-\ell+2}$ with coefficient
$$
\frac{k_1 \cdots k_\ell}{|\Aut\{k_1, \dots, k_\ell \}|}
$$
and then gets multiplied by $(m-\ell+3) \psi$, \dots, $m \psi$,
until it becomes
$$
\frac{k_1 \cdots k_\ell}{|\Aut\{k_1, \dots, k_\ell \}|}
\frac{m!}{(m-\ell+2)!} \psi^{\ell-2} i_{k_1, \dots, k_\ell}.
$$
Since the class $\psi^{\ell-2}$ is precisely the class
of a point in $\oM_{0;\ell+1}$ (see, for instance,~\cite{ELSV}),
we obtain precisely
the same coefficient as in the formulation of the theorem.

\paragraph{The expression for $\alpha_s \delta_{m_1, \dots, m_s}$.}

Recall that $\alpha_\ell \in H^{\ell-2}(\oM_{0;\ell+1})$ is the class
of a point.

According to Proposition~\ref{Prop:substitution},
to obtain the expression of $\alpha_s \delta_{m_1, \dots, m_s}$
we must substitute the expressions for $\psi^{m_i}$ into
the marked $\ell$-tree $T$:
\begin{center}
\
\includegraphics{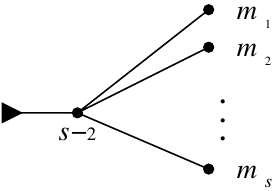}

\end{center}
It is easy to determine when a substitution of $\ell$ marked trees
$T_1, \dots, T_\ell$ into~$T$ contributes to the class
$\alpha_\ell i_{k_1, \dots, k_\ell}$. First of all, the marked trees $T_j$
must have $\ell$ leaves altogether with markings $k_1-1, \dots,
k_\ell-1$. Second, if a marked tree $T_j$ has $\ell_j$ leaves
with markings $k_1^{(j)}-1 , \dots, k_{n_j}^{(j)}-1$, then
it must describe a singularity class proportional to
$\alpha_{n_j} i_{k_1^{(j)}, \dots, k_{n_j}^{(j)}}$.
The coefficient of the singularity class
$\alpha_{n_j} i_{k_1^{(j)}, \dots, k_{n_j}^{(j)}}$
in the expression for $\nu_j^{m_j}$ is equal to the
coefficient of the monomial $x_{k_1^{(j)}} \cdots x_{k_{n_j}^{(j)}}$
in the polynomial $X_{m_j}$. Therefore, as claimed in the theorem,
the coefficient of the singularity class
$\alpha_\ell i_{k_1, \dots, k_\ell}$
in the expression for $\alpha_s \delta_{m_1, \dots, m_s}$ is equal to the
coefficient of the monomial $x_{k_1} \cdots x_{k_\ell}$ in the
polynomial $X_{m_1} \cdots X_{m_s}$.

\section{Completed cycles}
\label{Sec:CompCyc}

\subsection{Completed cycles and the classes $\psi^m$}

Let $\C S_N$ be the group algebra of the symmetric
group $S_N$. If $f \in \C S_N$ is a central element of the
group algebra, and $\lambda$ a partition of~$N$,
we denote by $f(\lambda)$ the scalar by which $f$
acts in the irreducible representation~$\lambda$.
(By Schur's lemma every central element of $\C S_N$
acts by a scalar in every irreducible representation of $S_N$.)
This identifies the center of $\C S_N$ with the algebra
of functions on the set of partitions of~$N$.
Let $C_{k_1, \dots, k_\ell}$ be the sum
of all permutations in $S_N$ with $\ell$ numbered cycles
of lengths $k_1, \dots, k_\ell$ respectively, and $N - \sum k_i$
non-numbered fixed points. For instance, $C_2$
is the sum of all transpositions, $C_{1,1}$ is $N(N-1)$
times the identity permutation, etc. Note that $C_{k_1, \dots, k_\ell}$
defines a central element in $\C S_N$ {\em simultaneously for all~$N$}.
For instance, the sum of transpositions $C_2$ is a well-defined
element of $\C S_N$ for any~$N$. If $N < \sum k_i$, then the corresponding
central element vanishes. Thus $C_{k_1, \dots, k_\ell}$
is actually a family of central elements in the group
algebras $\C S_N$ for all~$N$.

\begin{definition} \label{Def:StabCentElmt}
Given the positive integers $k_1, \dots, k_\ell$,
the family of central elements $C_{k_1, \dots, k_\ell} \in \C S_N$
for $N \geq 1$ is called a {\em stable\footnote{This has
nothing to do with stable maps.} central element}.
\end{definition}

\begin{theorem}[Kerov, Olshanski~\cite{KerOlsh}]
For every $m \geq 0$ there exists a unique linear
combination $\Cbar_{m+1}$ of stable
central elements $C_{k_1, \dots, k_\ell}$
such that
$$
\Cbar_{m+1}(\lambda) = \frac1{(m+1)!} \left[
(\lambda_i -i +1/2)^{m+1} - (-i+1/2)^{m+1}
\right]
$$
for every partition $\lambda = (\lambda_1 \geq \lambda_2 \geq \dots)$
of every integer~$N$.
\end{theorem}

\begin{example}
We have
\begin{eqnarray*}
\Cbar_1 &=& C_1,\\
\Cbar_2 &=& C_2,\\
\Cbar_3 &=& \frac12 C_3 + \frac14 C_{1,1} + \frac1{24} C_1,\\
\Cbar_4 &=& \frac16 C_4 + \frac13 C_{1,2} + \frac5{24} C_{2},\\
\Cbar_5 &=& \frac1{24} C_5 + \frac18 C_{1,3} +
\frac1{12} C_{2,2} + \frac1{36} C_{1,1,1}
+ \frac{11}{48} C_3+ \frac1{32} C_{1,1}+ \frac1{1920} C_1.
\end{eqnarray*}
\end{example}

\begin{definition}
The linear combination $\Cbar_{m+1}$ is called
the {\em completed $(m+1)$-cycle}.
\end{definition}

\begin{remark}
The set $S_\lambda = \{ (\lambda_i -i + 1/2)_{i\geq 1} \}
\subset \Z + \frac12$ is
a uniform way to encode a partition of any size. The empty
partition corresponds to the set $S_\emptyset$ of negative
half-integers. The set $S_\lambda$ differs from $S_\emptyset$
by a finite number of elements. The function $\Cbar_{m+1}$
is basically the renormalized sum of $(m+1)$-st powers
of the elements of $S_\lambda$ divided by $(m+1)!$.
\end{remark}

Denote by $S(z)$ the power series
$$
S(z) = \frac{\sinh(z/2)}{z/2}.
$$

\begin{proposition}[Okounkov, Pandharipande, \cite{OkoPan}, Proposition~3.2]
\label{Prop:CompCyc}
We have
$$
\Cbar_{m+1} =  \sum_{
\substack{g \geq 0, \\n \geq 1}
} \frac1{\ell!}
\sum_{
\substack{k_1, \dots, k_\ell \\
\sum k_i + \ell + 2g-2 = m}
}
\rho_{g;k_1, \dots, k_\ell}
C_{k_1, \dots, k_\ell},
$$
where the rational constant
$\rho_{g;k_1, \dots, k_\ell}$ is the coefficient of $z^{2g}$
in the power series
$$
\frac{\prod k_i}{K!}
S(z)^{K-1} \prod_{i=1}^\ell S(k_iz), \qquad K = \sum_{i=1}^\ell k_i.
$$
\end{proposition}

\begin{remark}
Our normalization of completed cycles $\Cbar_{m+1}$ differs
from that of~\cite{OkoPan} in three ways.
First, by a factor of $m!$; second, by the absence of a
constant term corresponding to $\ell=0$; third, by the
fact that, contrary to~\cite{OkoPan},
for us $C_{k_1, \dots, k_\ell}$ is the sum of
permutations with $\ell$ {\em numbered cycles}, which
changes the corresponding coefficient in the completed
cycle by a factor $|\Aut \{k_1, \dots, k_\ell \}|$.
\end{remark}

Every constant $\rho_{g;k_1, \dots, k_\ell}$ is associated to a
singularity of a higher genus stable map. Namely,
by analogy with
Definition~\ref{Def:I}, we can say that $f$ presents {\em
a singularity of type $I_{g;k_1, \dots, k_\ell}$} at $x \in C$
if $x$ lies on a contracted part of genus~$g$ meeting
$\ell$ branches of $C$ at ramification points of
orders $k_1, \dots, k_\ell$. It is well-known that if a stable
map $f$ presents a singularity of type $I_{g;k_1, \dots, k_\ell}$
at $x$, then the image $f(x)$ must be considered a
branch point of multiplicity $m=\sum k_i + \ell + 2g-2$.
In particular, if the stable map can be deformed into
a generic smooth map, the branch point will split
into $m$ simple branch points.

Since in this paper we only study stable maps of genus~$0$,
we will be only interested in the genus~$0$ part of the completed
cycles, that is, the part that corresponds to the $g=0$
terms in the sum of Proposition~\ref{Prop:CompCyc}, or in other
words, to the terms satisfying $\sum k_i + \ell = m+2$.

\begin{proposition}
\label{Prop:equality1}
If $\sum k_i + \ell = m+2$, then the coefficient
of $C_{k_1, \dots, k_\ell}$ in $\Cbar_{m+1}$ is the
same as the coefficient of $\alpha_\ell i_{k_1, \dots, k_\ell}$
in the expression of $\psi^m$ via the singularity classes,
where $\alpha_\ell$ is the class of a point in $\oM_{0;\ell+1}$.
\end{proposition}

\paragraph{Proof.} According to Theorem~\ref{Thm:CompCycles}
and Proposition~\ref{Prop:CompCyc} both are equal to
$$
\frac1{|\Aut\{ k_1, \dots, k_\ell\}|} \cdot
\frac{\prod k_i}{(\sum k_i)!}.
$$
\qed

\begin{remark}
In~\cite{OkoPan} Okounkov and Pandharipande established
a relation between Gromov-Witten invariants of $\CP^1$
and the completed cycles. Therefore the result of
Proposition~\ref{Prop:equality1} was to be expected.
Note, however, that Gromov-Witten invariants are only
intersection {\em numbers}, while here we are working
with {\em cohomology classes.} Therefore we can expect to get more
information than what is contained in the completed cycles.
And indeed, our expressions for $\psi^m$ involve
terms that do not appear in the completed cycles and
do not contribute to the computation of Gromov-Witten
invariants.
\end{remark}

\subsection{Products of completed cycles and the classes
$\alpha_s \delta_{m_1, \dots, m_s}$}

\begin{proposition}[\cite{KerOlsh}]
The product of two stable central elements is a finite
linear combination of stable central elements.
\end{proposition}

\begin{example}
We have $C_2^2 = C_{2,2} + 3 C_3 + \frac12 C_{1,1}$.
Indeed, $C_2$ is the sum of all transpositions, so
$C_2^2$ is the sum of products of all possible pairs
of transpositions. A product of two transpositions can
give either a permutation with two disjoint $2$-cycles (in a unique
way if we number the $2$-cycles in the same order as the
transpositions), or a $3$-cycle (and every $3$-cycle can
be decomposed into a product of two transpositions in three
possible ways), or an identity permutation with two marked
elements (but we can number these marked element in two
different ways).
\end{example}

The proof in the general case is a simple generalization
of the above example and is left to the reader.

\begin{proposition}
Let $\alpha_\ell \in H^{\rm top}(\oM_{0;\ell+1})$
be the class Poincar\'e dual to a point. Choose
$m_1, \dots, m_s$ and $k_1, \dots, k_\ell$ such that
$2s+ \sum m_j = \ell + \sum k_i$. Then the coefficient
of $\alpha_\ell i_{k_1, \dots, k_\ell}$ in the expression
of $\alpha_s \delta_{m_1, \dots, m_s}$ in terms of
singularity classes equals the coefficient of the
stable central element $C_{k_1, \dots, k_\ell}$ in the
product $\Cbar_{m_1+1} \cdots \Cbar_{m_s+1}$.
\end{proposition}

\paragraph{Proof.} Accroding to Theorem~\ref{Thm:CompCycles},
the coefficient of $\alpha_\ell i_{k_1, \dots, k_\ell}$ in the
expression for $\alpha_s \delta_{m_1, \dots, m_s}$ is equal to
the coefficient of the monomial $x_{k_1} \cdots x_{k_\ell}$ in
the product of polynomials $X_{m_1} \cdots X_{m_s}$. Recall
that the polynomials $X_m$ are defined as
$$
X_m = \sum_\ell \frac1{\ell!}
\sum_{
\substack{k_1, \dots, k_\ell\\
\sum k_i = m-\ell+2
}}
\frac{m!}{(m-\ell+2)!} \prod_{i=1}^\ell k_i x_i.
$$
They are transformed into the genus~$0$ part of
the completed cycles $\Cbar_m$ if we replace every monomial
$x_{k_1} \cdots x_{k_\ell}$ by the stable central element
$C_{k_1, \dots, k_\ell}$.

Let us call $\ell+ \sum k_i$ the {\em order} of the
stable central element $C_{k_1, \dots, k_\ell}$.
Then the genus~$g$ terms of a completed
cycle $\Cbar_m$ have order $m+2-2g$. The genus~$0$
elements have the biggest possible order $m+2$.

Consider two stable central elements $C_{k_1, \dots, k_\ell}$
and $C_{r_1, \dots, r_s}$. Denote their orders by $d_1$ and
$d_2$. Then we have
$$
C_{k_1, \dots, k_\ell} \cdot C_{r_1, \dots, r_s}
= C_{k_1, \dots, k_\ell,r_1, \dots, r_s} +
\mbox{ lower order terms},
$$
where the order of the first term equals $d_1+d_2$.
We conclude that {\em if we only keep the highest order
terms, then stable central elements multiply like monomials}:
$$
(x_{k_1} \cdots x_{k_\ell}) \;  \cdot \; (x_{r_1} \cdots x_{r_s})
\; = \; x_{k_1} \cdots x_{k_\ell} x_{r_1} \cdots x_{r_s}.
$$

This is enough to prove the second assertion of
Theorem~\ref{Thm:CompCycles}. Indeed, we have already
identified the highest order terms of a completed
cycle $\Cbar_{m+1}$ with the coefficients of the polynomial
$X_m$. Therefore the highest order terms of the product
$\Cbar_{m_1+1} \cdots \Cbar_{m_s+1}$ are identified with
the coefficients of the product of polynomials
$X_{m_1} \cdots X_{m_s}$.

\section{Appendix: sample computations}

\subsection{Expressions for $\prod (r \psi - \xi)$}
\label{Ssec:AppThm1}

Theorem~\ref{Thm:induction} leads to the following
expressions for $\prod (r \psi - \xi)$:

\begin{itemize}

\item
$\psi-\xi =$

$a_1$

\item
$(\psi-\xi)(2\psi-\xi)=$

$a_2 + \frac12 i_{1,1}$

\item
$(\psi-\xi)(2\psi-\xi)(3\psi-\xi)=$

$a_3 + 2 i_{1,2} + \frac16 i_{1,1,1}- \frac12 \xi i_{1,1}$.

\item
$(\psi-\xi)(2\psi-\xi)(3\psi-\xi)(4\psi-\xi)=$

$a_4 + 3 i_{1,3} + 2 i_{2,2} + i_{1,1,2} + \frac1{24} i_{1,1,1,1} +
\frac23 \psi i_{1,1,1}
-\xi (2 i_{1,2} + \frac16 i_{1,1,1})
+\frac12 \xi^2 i_{1,1}$

\item
$(\psi-\xi)(2\psi-\xi)(3\psi-\xi)(4\psi-\xi)(5\psi-\xi)=$

$a_5 + 4 i_{1,4} + 6 i_{2,3} + \frac32 i_{1,1,3} + 2 i_{1,2,2}
+\frac13 i_{1,1,1,2} + \frac1{120} i_{1,1,1,1,1} +
5 \psi i_{1,1,2}$

$+ \frac5{24} \psi i_{1,1,1,1}
-\xi(3 i_{1,3} + 2 i_{2,2} + i_{1,1,2} + \frac1{24} i_{1,1,1,1}
+ \frac32 \psi i_{1,1,1})
+\xi^2(2 i_{1,2} + \frac16 i_{1,1,1}) $

$- \frac12 \xi^3 i_{1,1}$

\end{itemize}

\subsection{Expressions for powers of $\psi$}
\label{Ssec:AppPsi}

By taking linear combinations of the
equalities of Section~\ref{Ssec:AppThm1},
we obtain the following expressions for the powers of $\psi$
(the underlined terms appear in Okounkov and Pandharipande's
completed cycles):

\begin{itemize}

\item
$\psi=\underline{a_1}+\xi$

\item
$\psi^2=\underline{\frac12 a_2} + \underline{\frac14 i_{1,1}}
+ \frac32 \xi a_1 + \xi^2$

\item
$\psi^3 = \underline{\frac16 a_3} + \underline{\frac13 i_{1,2}} +
\frac1{36}i_{1,1,1}+\frac{11}{12} \xi a_2+\frac38 \xi i_{1,1}
+ \frac74 \xi^2 a_1 + \xi^3$.

\item
$\psi^4 =
\underline{\frac1{24}a_4} +
\underline{\frac18 i_{1,3}} +
\underline{\frac1{12} i_{2,2}} + \frac1{24}i_{1,1,2} +
\frac1{576}i_{1,1,1,1} +
\underline{\frac1{36} \psi i_{1,1,1}}
+ \frac{25}{72} \xi a_3 + \frac{11}{18} \xi i_{1,2}$

$+ \frac{11}{216} \xi i_{1,1,1}
+\frac{85}{72} \xi^2 a_2 + \frac7{16} \xi^2 i_{1,1}
+\frac{15}8 \xi^3 a_1 + \xi^4$

\item
$\psi^5 = \underline{\frac1{120}a_5} +
\underline{\frac1{30} i_{1,4}} + \underline{\frac1{20} i_{2,3}}
+ \frac1{80} i_{1,1,3} + \frac1{60} i_{1,2,2}
+\frac1{360} i_{1,1,1,2} + \frac1{14400} i_{1,1,1,1,1}$

$+ \underline{\frac1{24} \psi i_{1,1,2}} +
\frac1{576} \psi i_{1,1,1,1}
+ \frac{137}{1440}\xi a_4 + \frac{25}{96} \xi i_{1,3}
+ \frac{25}{144} \xi i_{2,2} + \frac{25}{288} \xi i_{1,1,2}$

$+ \frac{25}{6912} \xi i_{1,1,1,1} + \frac{11}{216} \xi \psi i_{1,1,1}
+ \frac{415}{864} \xi^2 a_3 + \frac{85}{108} \xi^2 i_{1,2}
+ \frac{85}{1296} \xi^2 i_{1,1,1}
+ \frac{575}{432} \xi^3 a_2$

$+ \frac{15}{32} \xi^3 i_{1,1} + \frac{31}{16} \xi^4 a_1 + \xi^5$

\end{itemize}

\subsection{Basic via singularity classes}
\label{Ssec:AppBasViaSing}

Using the expressions for the powers of $\psi$ we obtain the
expressions for the other basic classes:

\paragraph{Codimension~2:}

\begin{itemize}

\item
$\delta_{0,0} = i_{1,1}$

\end{itemize}

\paragraph{Codimension~3:}

\begin{itemize}

\item $\delta_{0,0,0} = i_{1,1,1}$

\item $\delta_{0,1} = i_{1,2} + \xi i_{1,1}$

\end{itemize}

\paragraph{Codimension~4:}

\begin{itemize}

\item $\delta_{0,0,0,0} = i_{1,1,1,1}$

\item $\delta_{0,0,1} = i_{1,1,2} + \xi i_{1,1,1}$

\item $\psi \delta_{0,0,0} = \psi i_{1,1,1}$

\item $\delta_{1,1} = i_{2,2} + 2 \xi i_{1,2} + \xi^2 i_{1,1}$

\item $\delta_{0,2} = \frac12 i_{1,3} + \frac14 \psi i_{1,1,1}
+ \frac32 \xi i_{1,2} + \xi^2 i_{1,1}$

\end{itemize}

\paragraph{Codimension~5:}

\begin{itemize}

\item $\delta_{0,0,0,0,0} = i_{1,1,1,1,1}$

\item $\delta_{0,0,0,1} = i_{1,1,1,2} + \xi i_{1,1,1,1}$

\item $\alpha \delta_{0,0,0,0} = \alpha i_{1,1,1,1}$ for any class
$\alpha \in H^2(\oM_{0;5})$

\item $\delta_{0,1,1} = i_{1,2,2} + 2 \xi i_{1,1,2} + \xi^2 i_{1,1,1}$

\item $\delta_{0,0,2} = \frac12 i_{1,1,3} + \frac14
\left[
\begin{minipage}[c]{5em}
\includegraphics{delta012B.pdf}
\end{minipage}
\right]_\sing
+\frac32 \xi i_{1,1,2} + \xi^2 i_{1,1,1}$

\item $\psi \delta_{0,0,1} = \psi i_{1,1,2} + \xi \psi i_{1,1,1}$

\item $\delta_{1,2} = \frac12 i_{2,3} + \frac14 \psi i_{1,1,2}
+\frac12 \xi i_{1,3} + \frac32 \xi i_{2,2} + \frac14 \xi \psi i_{1,1,1}
+ \frac52 \xi^2 i_{1,2} + \xi^3 i_{1,1}$

\item $\delta_{0,3} = \frac16 i_{1,4} + \frac13 \psi i_{1,1,2} +
\frac1{36}
\left[
\begin{minipage}[c]{5em}
\includegraphics{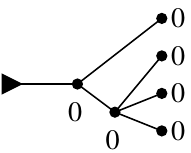}
\end{minipage}
\right]_\sing
+\frac{11}{12} \xi i_{1,3} + \frac38 \xi \psi i_{1,1,1}
+ \frac74 \xi^2 i_{1,2} + \xi^3 i_{1,1}
$

\end{itemize}

\subsection{Singularity via basic classes}
\label{Ssec:AppSingViaBas}

Inverting the relations of Sections~\ref{Ssec:AppPsi}
and~\ref{Ssec:AppBasViaSing} we obtain the following expressions
for the singularity classes in terms of basic classes:

\paragraph{Codimension~1:}

\begin{itemize}

\item $a_1 = \psi-\xi$

\end{itemize}

\paragraph{Codimension~2:}

\begin{itemize}

\item $i_{1,1} = \delta_{0,0}$

\item $a_2 = 2 \psi^2 -\frac12 \delta_{0,0} -3 \xi \psi + \xi^2$

\end{itemize}

\paragraph{Codimension~3:}

\begin{itemize}

\item $i_{1,1,1} = \delta_{0,0,0}$

\item $i_{1,2} = \delta_{0,1} - \xi \delta_{0,0}$

\item $a_3 = 6 \psi^3 - 2 \delta_{0,1} -\frac16 \delta_{0,0,0}
-11 \xi \psi^2 + \frac52 \xi \delta_{0,0} +6 \xi^2 \psi - \xi^3$

\end{itemize}

\paragraph{Codimension~4:}

\begin{itemize}

\item $i_{1,1,1,1} = \delta_{0,0,0,0}$

\item $i_{1,1,2} = \delta_{0,0,1} - \xi \delta_{0,0,0}$

\item $\psi i_{1,1,1} = \psi \delta_{0,0,0}$

\item $i_{2,2} = \delta_{1,1} - 2 \xi \delta_{0,1} + \xi^2 \delta_{0,0}$

\item $i_{1,3} = 2 \delta_{0,2} - \frac12 \psi \delta_{0,0,0}
- 3 \xi \delta_{0,1} + \xi^2 \delta_{0,0}$

\item $a_4 = 24 \psi^4 - 6 \delta_{0,2} - \delta_{1,1} -\delta_{0,0,1}
+\frac{13}6 \psi \delta_{0,0,0}- \frac1{24} \delta_{0,0,0,0}
-50 \xi \psi^3 + \frac15 \xi \delta_{0,1} +$

$\frac76 \xi \delta_{0,0,0}
+35 \xi^2 \psi^2 - \frac{15}2 \xi^2 \delta_{0,0} - 10\xi^3 \psi +\xi^4$

\end{itemize}


\begin{thebibliography}{99}

\bibitem{ELSV} {\bf T. Ekedahl, S. Lando, M. Shapiro, A. Vainshtein.}
{\em Hurwitz numbers and intersections on moduli spaces of curves.}
-- Invent. Math.~{\bf 146} (2001), no.~2, 297--327.

\bibitem{KazLan0} {\bf M. Kazaryan, S. Lando.} {\em On intersection theory on Hurwitz spaces.} (Russian) Izv. Ross. Akad. Nauk Ser. Mat. 68 (2004), no. 5, 91--122; translation in Izv. Math. 68 (2004), no. 5, 935--964.

\bibitem{KazLan} {\bf M. Kazarian, S. Lando.} {\em Thom
polynomials for maps of curves with isolated singularities.}
-- Tr. Mat. Inst. Steklova, {\bf 258} (2007), Anal. i Osob. Ch. 1, 93--106 (Russian) ; translation in Proc. Steklov Inst. Math., {\bf 258} (2007), no.~1,
\texttt{arXiv: 0706:1523}.

\bibitem{KLZ2} {\bf M. Kazarian, S. Lando, D. Zvonkine.} In preparation.

\bibitem{Keel} {\bf S. Keel.} {\em Intersection theory of moduli space of stable n-pointed curves of genus zero.}
-- Trans. Amer. Math. Soc., {\bf 330} (1992), no.~2, 545--574.

\bibitem{KerOlsh} {\bf S. Kerov, G. Olshanski.} {\em Polynomial functions on the set of Young diagrams.}
-- C. R. Acad. Sci. Paris, S\'er.~I, Math. {\bf 319} (1994), no.~2, 121--126.

\bibitem{OkoPan} {\bf A. Okounkov, R. Pandharipande.}
{\em Gromov-Witten theory, Hurwitz theory, and completed cycles.}
Ann. of Math.{\bf 163} (2006), no.~2, 517--560.

\bibitem{Thom} {\bf R. Thom.} {\em Quelques propri\'et\'es
globales des vari\'et\'es diff\'erentiables.}
-- Comment. Math. Helv.,~{\bf 28} (1954), 17--86.

\end{thebibliography}
\end{document}